\documentclass[10pt,a4paper]{article}
\usepackage{latexsym,amssymb}
\usepackage{enumerate, verbatim, mathrsfs}
\usepackage{graphics}

\usepackage{amsthm}
\usepackage{amssymb}
\usepackage{latexsym}
\usepackage{longtable}
\usepackage{epsfig}
\usepackage{amsmath}
\usepackage{hhline}

\usepackage{latexsym,amssymb,verbatim}
\usepackage{enumerate}
\usepackage{graphics}

\newtheorem{thm}{Theorem}[section]

\newcommand{\pa}{{\partial}}

\newcommand{\Ker}{\mbox{{\rm Ker}\,}}

\newcommand{\dne}{\hfill $\Box$} 

\newcommand{\FF}{{\mathbb F}}
\newcommand{\CP}{{\cal P}}
\newcommand{\CQ}{{\cal Q}}
\newcommand{\spa}{\!\!\!\!\!\phantom{|^{|^{|}}_{|_{|^{|}}}}\!\!\!\!\!}

\begin{document}
%
%
\title{\vspace*{-15mm}
{\Huge\sf On Stanley's Inequalities for  \\ Character Multiplicities} \\}
\author{\\V B Mnukhin%
{\small\thanks{Support from the London Mathematical Society
and the Leverhulme Trust is acknowledged.}}\\
\\
{\normalsize Department of Mathematics and Computer Science,}\\
{\normalsize Southern Federal University Taganrog, Russia}\\
{\normalsize e-mail:{\sl \,\,mnukhin.valeriy@mail.ru}}\\[4mm]
\and
I J Siemons\\
{\normalsize School of Mathematics, University of East Anglia}\\
{\normalsize Norwich NR4 7TJ, United Kingdom.}
{\normalsize e-mail:{\sl \,\,J.Siemons@uea.ac.uk}}}
\date{\scriptsize Version of 1 August, 2011, printed \today}


\maketitle
%
\begin{abstract}\noindent
Let  $G$ be  a group of automorphisms  of a ranked poset $\CQ$ and 
let $N_{k}$ denote  the  number of orbits on the elements of rank $k$ in $\CQ.$ 
What can be said about the $N_{k}$ for standard posets, such as finite projective spaces or the Boolean lattice? We discuss the connection of this question to the representation theory of the group, and in particular to the inequalities of Livingstone-Wagner and Stanley.
We show that these are special cases 
of more general inequalities which depend on the prime divisors of the group order. The new inequalities often yield stronger bounds depending on the  order of the group.

\bigskip\noindent {\sc MSC}:\quad Finite permutation groups \,20B05,
Automorphism groups of combinatorial structures \,20B25,
Group actions on posets and homology groups of posets  \,05E25.

\end{abstract}
\maketitle

\section{\sc Introduction}                %

The enumeration of combinatorial structures up to a group action 
belongs to the oldest problems in combinatorics. Often this means 
that one has to estimate the number of orbits of a permutation group 
acting on a poset, and more specifically, on the lattice of subsets of a set. 
If $G$ is a permutation group on the set $\Omega$ of size $n$ 
let therefore $N_k$ be the number of orbits of $G,$ in the natural action on the 
$k$-subsets of $\Omega$. What can be said about these numbers? 
One answer is given by P\'olya's enumeration theory, namely that 
$$
\sum_{k=0}^{\infty} N_{k}t^k={\rm Z}_G(1+t,1+t^2,\ldots,1+t^n)
$$
when ${\rm Z}_G$ is the {\sl cycle index} of the group.
However, usually ${\rm Z}_G$ is quite complicated and so
this formula is often  inefficient. 

\medskip
From elementary group theory another answer is on offer: $N_{k}$ is the average number of $k$-element subsets of $\Omega$  fixed by $G$. This means that
$$N_{k}=\frac{1}{|G|}\,\sum_{g\in G}\,\, {\rm fix}_{k}(g)$$ where
${\rm fix}_{k}(g)$ is the number of $k$-element subsets  that are fixed by the element $g\in G$.
Also this formula is difficult to evaluate in general. However, it has a convenient interpretation as  the inner product of the permutation character  ${\rm fix}_{k}(g)$ with the identity character of $G.$ Thus if we work with characters over $\mathbb{C}$ then $N_{k}$ is the multiplicity of the identity character in the permutation character of $G$ acting on $k$-element subsets.

As the inclusion matrix between $l$- and $k$-subsets of $\Omega$  over the rationals  has rank ${n\choose k}$ for all $l\geq k$ and $k+l\leq n$, it follows from general representation 
theory that 
$$N_{l}\geq N_{k}\,\,\mbox{\quad if $l\geq k$ and $l+k\leq n$.}\eqno (1)$$
This inequality is due to Livingstone and Wagner. 
In~\cite{Stanley} Stanley \,observes that
$$c_{l}\geq c_{k}\,\,\mbox{\quad if $l\geq k$ and $l+k\leq n$}\eqno (2)$$
when $c_{k}$ denotes the multiplicity of a given fixed complex irreducible character of $G$ in the permutation character on $k$-element subsets, see also~\cite[Section 9.4]{Dixon} and~\cite[Chapter 9.4]{Sagan}. 

Suppose more generally that $G$ is an automorphism groups of a ranked partially ordered set $\CQ.$ For $0\leq k\in \mathbb{Z}$  let ${\cal Q}_{k}$ denote the set of elements of rank $k$ in $\CQ.$   
Then it is clear that  \,(2)\, holds for the action of $G$ on any two rank sets  ${\cal Q}_{k}$ and ${\cal Q}_{l}$ provided that  the incidence matrix between ${\cal Q}_{k}$ and ${\cal Q}_{l}$ has rank equal to $|{\cal Q}_{k}|$  over the rationals.  This applies for instance to groups preserving a simplicial complex or to  the action of a linear group on the Grassmanians of the vector space.

\medskip

In this paper we establish  corresponding characters inequalities in positive characteristic. Let $F$ be a field of characteristic $p>0.$ We suppose that $G$ is either a permutation group on an $n$-set $\Omega,$ or  a linear group  $G\subseteq {\rm P}\Gamma{\rm L}(n,q)$ with $q$ not divisible by $p.$ Put $\pi=p$ or $\pi=\pi(p,q),$ respectively, where $\pi(p,q)$ is the characteristic of $q$ in $F,$ see the definition in Section 3. If $k\mapsto c_{k}$  is an integer  function  put $[c_{k}]_{\pi}:=\sum_{t\in\mathbb{Z} }\,c_{k+t\pi}$ whenever this is defined.  The main result is  the following

\medskip
\begin{thm} \, Let $G$ and $F$ be as given, and suppose that $U$ is an irreducible $FG$-module. Let $c_{k}$ denote the multiplicity of $U$ in the permutation representation  of $G$ on the $k$-element subsets of $\Omega,$ or the $k$-dimensional subspaces of ${\rm GF}_{q}^{n},$ respectively. Then 
$$
\bigl[ c_{m}\bigr]_\pi\ge
\bigl[ c_{m-1}\bigr]_\pi\ge
\bigl[ c_{m-2}\bigr]_\pi\ge
\ldots\ge
\bigl[ c_{m-s}\bigr]_\pi\ge 0
$$
where $m=\lfloor\frac{n}{2}\rfloor$ and $s=\lfloor\frac{\pi}{2}\rfloor$.

\end{thm}
We have included some applications in Section 4. Here just note that the inequality in the theorem becomes Stanley's inequality \,(2)\, when $\pi>n.$ Another application concerns the orbit numbers $N_{k}$ when $p$ does not divide the group order. In this case  the theorem implies that $$\bigl[ N_{m}\bigr]_\pi\ge
\bigl[ N_{m-1}\bigr]_\pi\ge
\bigl[ N_{m-2}\bigr]_\pi\ge
\ldots\ge
\bigl[ N_{m-s}\bigr]_\pi\ge 0.
$$ This is an improvement on the Livingstone--Wagner inequality \,(1)\, when information about the group order is available. Examples show that the new orbit inequalities give bounds which are sharp in some cases. The theorem has also implications  for multiplicity-free representations in positive characteristics.

\section{\sf Homology and the Hopf-Lefschetz Formula}%

\medskip
Let $F$ be a field and for each $k\in \mathbb{Z}$ let $M_k$ be a finite dimensional vector space over $F$. We assume that there is a positive integer $n$ so that $M_{k}=0$ if $k>n$ or $k<0,$ and we 
 set $M=\bigoplus_{k\in \mathbb{Z}}\,M_{k}.$  If $\partial: \,\, M\rightarrow  M$ is any linear map with $\partial(M_{k})\subseteq M_{k-1}$ for all $k$ then we may investigate the sequence
$$
{\cal M }:\quad
 0\stackrel{\,\partial}{\longleftarrow} M_0\stackrel{\,\partial}{\longleftarrow}
 M_1\stackrel{\,\partial}{\longleftarrow}
\ldots\stackrel{\,\partial}{\longleftarrow}
M_{j}\stackrel{\,\partial}{\longleftarrow}
\ldots\stackrel{\,\partial}{\longleftarrow}
 M_n \stackrel{\,\partial}{\longleftarrow}
 0\longleftarrow
 \cdots\,\,.
$$
Note that $\partial$ is nilpotent since $\pa^{n+1}=0.$  Let therefore  $\pi>1$ be the least integer such that $\partial^{\pi}=0.$ When $\pi=2$ then ${\cal M }$ is a homological sequence in the usual 
sense. 

When $\pi>2$ then the  methods of classical homology theory can be adjusted by focusing on certain sequences obtained from 
${\cal M},$ as follows.  For any $j$ and  $0<i<\pi$ consider the sequence $$
{\cal M}^{}_{(l,r)}:\quad \ldots\stackrel{\,\,\,\partial^*}{\longleftarrow}
M_{j-\pi}^{}\stackrel{\,\,\,\partial^*}{\longleftarrow}
 M_{j-i}^{}\stackrel{\,\,\,\partial^*}{\longleftarrow}
M_{j}^{}
\stackrel{\,\,\,\partial^*}{\longleftarrow} M^{}_{j+\pi-i}
\stackrel{\,\,\,\partial^*}{\longleftarrow} M^{}_{j+\pi}
\stackrel{\,\,\,\partial^*}{\longleftarrow}\ldots\eqno(3)
$$
in which $\partial^{*}$ is the appropriate power of $\partial$. Since $\partial^{\pi}=0$ we have $\partial^{*}\partial^{*}=0$ and therefore ${\cal M}^{}_{(l,r)}$ is a homological sequence. We denote the
homology at $M_{j-i}^{}\leftarrow M_{j}^{}\leftarrow M^{}_{j+\pi-i}$
by
$$
H_{j,i}:=(\Ker\partial^{i}\cap M_{j}^{})/
                      \partial^{\pi-i}(M_{j+\pi-i}^{})\,\,.
$$
We may say that  ${\cal M}$ is {\it $\pi$-homological} and that $H_{j,i}$ is  the {\it generalized homology module} of ${\cal M}$ for parameters $j$ and $0<i<\pi.$ For $\pi=2$ these are just the customary definitions.  

In order to keep track of the parameters in a sequence we need the following definitions. The
sequence \,(3)\, is determined by any arrow $M^{}_{l}\leftarrow M^{}_{r}$ in
it and hence it is denoted by ${\cal M}^{}_{(l,r)}.$ (Of course, $(l,r)$
stands for left-right.) The unique arrow with $0\le a+b<\pi$ is the 
{\it initial arrow\/}. The module $M^{}_{b}$ in the initial arrow
is the $0$-{\it position\/} of ${\cal M}^{}_{(l,r)}.$ The position of any 
other module in ${\cal M}^{}_{(l,r)}$ is counted from this $0$-position, 
and $(a,b)$ is the {\it type\/} of ${\cal M}^{}_{(l,r)}.$

\bigskip
We now state an appropriate  version of the Hopf-Lefschetz  trace formula for ${\cal M}^{}_{(l,r)}$ in the presence of  a group action.  Suppose that the group $G$  acts linearly on each of the $M_{k}$ and that $\partial$ commutes with  this action. In this case ${\cal M}$ is an  $FG$-sequence and all $H_{j,i}$ are  $FG$-modules. 

When ${\cal M}$ is homological, that is $\partial^{2}=0,$ then the  standard Hopf-Lefschetz formula says that the virtual modules $\sum_{k} \,\, (-1)^{k}M_k$ and $\sum _{k}\,\,(-1)^{k}H_k$ are the same, see for instance~\cite[Chapter 2,\,\S22]{Munkres}. More generally, if  $A$ and $B$ are virtual $FG$-modules (that is, formal combinations of $FG$-modules with coefficients in $\mathbb Z),$ then $A$ is {\it equivalent to $B$,} denoted $A\sim B$, if
all irreducible $FG$-modules appearing in $A$ and $B$ have the same integer multiplicity. Therefore the Hopf-Lefschetz formula takes the form  $\sum_{k} \,\, (-1)^{k}M_k\,\sim\,\sum _{k}\,\,(-1)^{k}H_k.$
Adjusting this carefully to the terms in  \,(3)\, we have

\medskip
\begin{thm}\label{HopfLefschetz}
\,{\rm (Hopf-Lefschetz Trace Formula):}\quad Let $G$ be a group and let 
${\cal M}$ be a $\pi$-homological $FG$-sequence of finite-dimensional 
modules. Then 
$$
\sum_{t\in\mathbb{Z}} \,\,\, (M_{j-\pi t} -M_{j-i-\pi t})
\,\,\sim\,\,\sum_{t\in\mathbb{Z}}\,\,\,(H_{j-\pi t} -H_{j-i-\pi t})
$$
as virtual $FG$-modules,  for any $j$ and $0<i<\pi$. 
\end{thm}

\medskip
This formula is particularly useful when many of the homologies vanish. We say that  ${\cal M}^{}_{(l,r)}$ is 
{\it almost exact } if \,(3)\, has at most one non-vanishing homology, which then is 
denoted by $H_{(l,r)}.$ If this homology appears in position $d$ then the trace formula becomes 
$$H_{(l,r)}\,\,\sim\,\,(-1)^{d}
\sum_{t\in\mathbb{Z}} \,\,\,\,\, (M_{b-\pi t} -M_{a-\pi t}) $$
where $(a,b)$ is the type of ${\cal M}^{}_{(l,r)}.$ To simplify such sums, if $f$ is a function on $\mathbb{Z}$   we abbreviate  $[f_{j}]_{\pi}:=\sum_{t\in\mathbb{Z}} \,\,\, f(j-\pi t)$ and extend this  linearly by setting   $[f_{j}-f_{k}]_{\pi}=[f_{j}]_{\pi}-[f_{k}]_{\pi}$.  Our results will be an application of the following version of the trace formula.

\medskip
\begin{thm}\label{HopfLefschetz2}
\quad
Let $G$ and ${\cal M}$ be as above and suppose that ${\cal M}^{}_{(l,r)}$ is almost exact for some $(l,r)$ with non-trivial homology in position $d.$ Then $$H_{(l,r)}\,\,\sim\,\,(-1)^{d}\,\,[M_{b}-M_{a}]_{\pi}$$ as  virtual $FG$-modules where $(a,b)$ is the type of ${\cal M}^{}_{(l,r)}$.
\end{thm}

\section{\sf Incidence Homology in some Standard Posets}
The situation just described occurs commonly  in combinatorics. Let $\cal Q=(\cal Q,\leq)$ be a finite ranked poset. For $k\in \mathbb{Z}$ we denote the set of all $x\in {\cal Q}$ with ${\rm rk}(x)=k$ by   ${\cal Q}_{k}.$ If  $F$ is a field of characteristic $p\geq 0$ let $M_{k}=F{\cal Q}_{k}$ be the vector space with basis ${\cal Q}_{k}.$  Consider the linear {\it incidence map} $\partial\!: M_{k}\to M_{k-1}$ which is defined for $x\in  {\cal Q}_{k}$ by $$\partial(x)=\sum y$$ where the sum runs over all $y\in{\cal Q}_{k-1}$ with $y\leq x.$ Evidently, if $G$ is a group of automorphisms  of ${\cal Q}$ then $\partial$ is an $FG$-homomorphism on $M=\bigoplus_{k\in\mathbb{Z}}\,M_{k}.$ Note that $\partial$ is nilpotent as $\partial^{m}=0$ for $m=1+\max\{{\rm rk}(x)\,\,:x\in {\cal Q}\}< \infty$ since $\cal Q$ is finite.

The homology modules $H_{j,i}$ defined above for $(M,\partial)$ will be called the {\it incidence homologies} of $\cal Q$ over $F.$ The most interesting examples of homology occur when $F$ has characteristic $p>0.$ Here  the situation is well understood  if ${\cal Q}$ is a standard poset such as the power set of a finite set or a finite projective space. To extend this to results on wider classes of posets would certainly be interesting. 

\vspace{0.3cm}
{\large\sc Projective Spaces:}\quad Let $q$ be a prime power {\it not} divisible by $p$ and let $V=\FF^{n}$ where  $\FF:={\rm GF}_{q}.$  We denote the projective space associated to $V$ by $\CP(n,q).$ So $\CP(n,q)$  is the poset of all subspaces of $V$ ordered by inclusion; it is ranked by the dimension of subspaces and has automorphism group $G_{n}:={\rm P}\Gamma{\rm L}(n,q).$

Now apply the definitions above to  $\CQ=\CP(n,q).$  For the positive integer $i$ let  $
|i|_q:=1+q+q^2+...+q^{i-1} $ and $(i!)_q :=|1|_q\cdot|2|_q\cdots|i|_q\,\,.$ It is easy to see that for $x\in M_{k}$ we have $$\pa^{i}(x)= (i!)_q \sum\,\,y$$ where the sum runs over all $y\in \CQ_{k-i}.$ (The coefficient $(i!)_q$ counts the number of saturated chains $y=y_{i}<...<y_{1}<y_{0}=x$ for a fixed $y\in \CQ_{k-i}.$) Now we are interested in the least integer $i=\pi(p,q)$ so that $(i!)_q=0$ in $F.$

\medskip
{\sc Definition} \quad
{\it For the prime $p$ and the  integer $q$  co-prime  to $p$ the {\rm characteristic}  of $q$ in ${\rm GF}_{p}$ is the least integer  $\pi=\pi(p,q)>0$ so that  $|\pi|_q \equiv 0$ mod $p.$}

\medskip
For $\pi(p,q)$ also the name {\it quantum characteristic} is in use. If $q\equiv 1$ mod $p$ then $\pi(p,q)=p$ while if $p$ does not
divide $q-1$ then $\pi$ is the order of $q$ modulo $p,$ since
$(1+q+q^{2}+...+q^{\pi-1})=(q^\pi -1)/(q-1)$. In either case
$q^\pi\equiv 1$ mod $p$ and if
$\pi\geq 2$ then $((\pi-1)!)_q\not\equiv 0$ mod $p$ while
$(\pi!)_q\equiv 0$ mod $p$. The Table~1 reproduced from \cite{Projective} gives the  value of $\pi(p,q)$ for small  $p$ and $q$ as an  illustration. 

{\footnotesize
\begin{center}
\begin{tabular}{|c||r|r|r|r|r|r|r|r|r|r|r|r|r|}
\hline
{\normalsize $\,\,\,\spa \pi(p,q)\spa\,\,\,$} & \multicolumn{13}{|c|}{Values for $q$} \\
\cline{1-14}
$p$& 2  & 3  & 4  & 5  & 7  & 8  & 9  & 11 & 13
   & 16 & 17 & 19 & \,23\,  \\
\hline
2  &-- & 2  &-- & 2  & 2  &-- & 2  & 2  & 2
   &-- & 2  & 2  & 2   \\
3  & 2  &-- & 3  & 2  & 3  & 2  &-- & 2  & 3
   & 3  & 2  & 3  & 2   \\
5  & 4  & 4  & 2  &-- & 4  & 4  & 2  & 5  & 4
   & 5  & 4  & 2  & 4   \\
7  & 3  & 6  & 3  & 6  &-- & 7  & 3  & 3  & 2
   & 3  & 6  & 6  & 3  \\
11 & 10 & 5  & 5  & 5  & 10 & 10 & 5  &-- & 10
   & 5  & 10 & 10 & 11 \\
13 & 12 & 3  & 6  & 4  & 12 & 4  & 3  & 12 &--
   & 3  & 6  & 12 & 6   \\
17 & 8  & 16 & 4  & 16 & 16 & 8  & 8  & 16 & 4
   & 2  &-- & 8  & 16  \\
19 & 18 & 18 & 9  & 9  & 3  & 6  & 9  & 3  & 18
   & 9  & 9  &-- & 9  \\
\hline
\end{tabular}
\vspace{2mm}
\end{center}}
\centerline{\sc Table 1:  The Function $\pi(p,q)$}

In particular, if we let $j$ be arbitrary and $0<i<\pi(p,q)$ then the sequence ${\cal M}_{(l,r)}$ in \,(3)\, above  is homological. The next result has been proven in~\cite{Projective}\,:

\medskip
\begin{thm}\label{proj}\quad
Let $F$ be a field of characteristic $p>0$ where $p$ does not divide $q.$ Then for all $j\geq 0$ and all $i$ with $0<i<\pi(p,q)$ we have $H_{j,i}=0$
unless $n-\pi<2j-i<n.$ In particular, for fixed $j$ and $i$ the sequence ${\cal M}_{(l,r)}$ in \,${\rm (3)}$\,  is almost exact. 
\end{thm}

The exact $FG_{n}$-structure of the homology modules has been determined in \cite{HomGrass}, this is however not important for the current paper.

\vspace{0.3cm}
{\large\sc Power Sets:}\quad We extend the definition above  and let $\CP(n,1)$ denote the lattice of all subsets of an $n$-set. This poset is ranked by the size of subsets and has automorphism group ${\rm Sym}_{n}.$ If we apply the definitions from the beginning of this section  to $\CQ=\CP(n,1)$ then we find that $\pa^{i}=0$ if and only if $i\geq p.$ Thus the least integer $\pi$ with $\pa^{\pi}=0$ is  
$\pi=\pi(p,1)=p$ according to the definition. In fact, for the purpose of the incidence homology of  $\CP(n,1)$ over ${\rm GF}_{p}$ we may regard  $\CP(n,1)$ as the projective space $\CP(n,q)$  over the `field' with $q=1$ element.

Also here, if we let $j$ be arbitrary and $0<i<p$ then the sequence ${\cal M}_{(l,r)}$  in \,(3)\, is homological. The precise equivalent to Theorem~\ref{proj}  \, has been proved in~\cite{Finite}:

\medskip
\begin{thm}\label{bool}\quad
Let $F$ be a field of characteristic $p>0.$ Then for all $j\geq 0$ and all $i$ with $0<i<p$ we have $H_{j,i}=0$ unless $n-p<2j-i<n.$ In particular, for fixed $j$ and $i$ the sequence ${\cal M}_{(l,r)}$ in \,${\rm (3)}$\, is almost exact. 
\end{thm}

The $F{\rm Sym_{n}}$-structure of the homology modules is completely analogue to the $q>1$ case mentioned earlier.

\section{\sc Group Actions}
We return to the situation when $G$ is a group of automorphisms of $\CQ$  where $\CQ=\CP(n,q)$ when $q$ is a prime power or $\CQ=\CP(n,1)$ when $\CQ$ is the power set of an $n$-set. Then $G$ preserves the rank sets $\CQ_{k}$ for all $k.$ Let $F$ be a field of characteristic $p$ not dividing $q,$ and correspondingly put $\pi=\pi(p,q)$ or $\pi=p$ when $q=1.$ Define $M_{k}:=F\CQ_{k}$ as before. In both cases it is a  key fact  that any sequence  of the kind ${\cal M}(l,r)$  as in \,(3)\, is  almost exact. Hence the trace formula gives useful information about character multiplicities in the $G$-representation on $\CQ_{k}$ in characteristic $p>0.$

\vspace{0.3cm}
{\large\sc Character Multiplicities:}\quad We fix some  $j\geq 0$ and $0<i<\pi.$ Then the  sequence 
$${\cal M}_{(l,r)}:\quad
\ldots{\leftarrow}
M_{j-\pi}^{}{\leftarrow}
 M_{j-i}^{}{\leftarrow}
M_{j}^{}
{\leftarrow} M^{}_{j+\pi-i}
{\leftarrow} M^{}_{j+\pi}
{\leftarrow}\ldots\eqno(4)
$$ is almost exact according to the theorems in Sections~2 and~3. Let its type be  $(a,b),$ with non-trivial homology $H_{(l,r)}$ in position $d$. If $G\subseteq {\rm Aut}(\CQ)$ let $U$ be an irreducible $FG$-module. For an arbitrary $FG$-module $M$ let $0\leq c(M,U)$ be the multiplicity of $U$ in $M.$ 

\medskip
\begin{thm}\label{Ineq1}
\,\, Let $q\geq 1$ and let $F$ be  field of characteristic $p>0$ with $p$ co-prime to $q.$ Suppose that $G\subseteq {\rm Aut}(\CP(n,q))$ and that $U$ is an  irreducible irreducible $FG$-module. For each module $M_{k}$  in~${\rm (4)}$ let $c_{k}:=c(M_{k},U).$ Then  $c_{k}=c_{n-k}$  and $$c(H_{(l,r)},U)\,=\,(-1)^{d}\,\,[c_{b}-c_{a}]_{\pi}\,\,.$$ In particular, \,$[c_{a}]_{\pi}\,\leq\, [c_{b}]_{\pi}$\, if $d$ is even and \,$[c_{b}]_{\pi}\,\leq\, [c_{a}]_{\pi}$\, if $d$ is odd.
\end{thm}

Note, if $\pi>n$ then $[c_{a}]_{\pi}=c_{a}$ and $[c_{b}]_{\pi}=c_{b}.$ Therefore, for $a\leq b$ and $a+b\leq n$ the second part of the theorem gives $c_{a}\leq c_{b}.$ This is  the exact form of  Stanley's inequality (2).

\medskip
{\it Proof:} \, Since $p$ and $q$ are co-prime $M_{k}$ and $M_{n-k}$ are isomorphic as $FG$-modules, by Theorem~14.3\, of James~\cite{James}. In particular,  $c_{k}=c_{n-k}.$ The second part   follows from the trace formula Theorem~\ref{HopfLefschetz2}\, and Theorems~\ref{proj}~and~\ref{bool}. The inequalities come from the fact that $0\leq c(H_{(l,r)},U).$\dne

\bigskip
For the full sequence $$
{\cal M }:\quad
 0\leftarrow M_0\leftarrow
 M_1\leftarrow
\ldots\leftarrow
M_{j}\leftarrow\ldots\leftarrow
 M_n \leftarrow
 0\eqno (5)
$$
we have ${\pi\choose 2}$ choices for the $(l,r)$-term in  a sequence ${\cal M}_{(l,r)}$ as in   \,(4), and hence  there are ${\pi\choose 2}$ inequalities as in Theorem~\ref{Ineq1}.  The resulting  inequalities can be summarized succinctly as follows.

\medskip

\begin{thm}\label{Ineq2} 
\,\, Let $q\geq 1$ and let $F$ be  field of characteristic $p>0$ with $p$ co-prime to $q.$ Suppose that $G\subseteq {\rm Aut}(\CP(n,q))$ and that $U$ is an  irreducible irreducible $FG$-module. For each module $M_{k}$  in~${\rm (4)}$ let $c_{k}:=c(M_{k},U).$ Then  $c_{k}=c_{n-k}$ and 
$$
\bigl[ c_{m}\bigr]_\pi\ge
\bigl[ c_{m-1}\bigr]_\pi\ge
\bigl[ c_{m-2}\bigr]_\pi\ge
\ldots\ge
\bigl[ c_{m-s}\bigr]_\pi\ge 0
$$
where $m=\lfloor\frac{n}{2}\rfloor$ (that is, either $n=2m$ or $n=2m+1$),
and $s=\lfloor\frac{\pi}{2}\rfloor$.
\end{thm}

\medskip
{\it Proof.\quad} When $\pi>n$ then  $[c_{k}]_\pi=c_{k}$ and we have the inequalities discussed above.  So we suppose that $\pi\le n$. From Theorem~\ref{Ineq1}\, we have   
$c_{k}=c_{n-k}$  and hence  $[c_{k}]_{\pi}=[c_{n-k}]_{\pi}$ for all $k$. 

First let $n=2m$. Now apply Theorem~\ref{Ineq1} \, for the pairs 
$(l,r)=(m-\pi+1,m),\, (m-1-\pi+3,m-1),\, (m-2-\pi+5,m-2),...$ \, resulting in the inequalities 
$$
\begin{array}{lcl}
\ [c_m]_\pi      & \ge &     [c_{m-\pi+1}]_\pi=[c_{m+1}]_\pi=[c_{m-1}]_\pi\; ,\\\
  [c_{m-1}]_\pi  & \ge &     [c_{m-1-\pi+3}]_\pi=[c_{m+2}]_\pi=[c_{m-2}]_\pi\;,   \\\
  \cdots         & \cdots & \cdots                              \\\
  [c_{m-s+1}]_p  & \ge &     [c_{m-s-1}]_p=[c_{m+s}]_p    \;,   \\\
[c_{m+s}]_p    & \ge &     [c_{m+s-p+1}]_p=[c_{m-s}]_p  \;
\end{array}
$$

\vspace{-0.2cm}
which prove the theorem for $n$ even. Secondly, if $n=2m+1$ choose the pairs 
$(l,r)=(m-\pi+2,m),\, (m-1-\pi+4,m-1),\, (m-2-\pi+6,m-2),...$ when the results follow in the same way. It is easy to note that these inequalities together with $c_{k}=c_{n-k}$ are the only significant inequalities, and other inequalities may be derived from the ones in the theorem. \dne

As an application  we show that the theorem can be used to prove that certain representations are not multiplicity free. Surprisingly this property depends on $n$ and $\pi$ only, without reference to the group. For instance, take $\CQ=\CP(10,2),$ let $G\subseteq {\rm P}\Gamma{\rm L}(10,2), $ let $F={\rm GF}_{p}$ with $\pi(p,2)=8$ and let $U$ be an irreducible $FG$-module. (From Table 1 we may take $p=17$ since $\pi(17,2)=8.)$ If $c_{k}$ is the multiplicity of $U$ in $M_{k}$ as in the theorem then $$[c_{5}]=c_{5}\geq [c_{4}]=c_{4}\geq [c_{3}]=c_{3}\geq [c_{2}]=c_{10}+c_{2}\geq[c_{1}]=c_{9}+c_{1}=2c_{1}$$ since $c_{9}=c_{1}.$ Hence, any irreducible $U$ that appears in the action of $G$ on $\CQ_{1}$ with multiplicity $c_{1}>0$ has multiplicity $\geq 2c_{1}$ in the action on $\CQ_{k}$
 for $k\in \{3,4,5,6,7\}.$ Furthermore, if $U\neq 1$ then $c_{10}=c_{0}=0$ and hence $U$ has also multiplicity $\geq 2c_{1}$ in the action on $\CQ_{2}$ and $\CQ_{8}.$ It would be easy but tedious to formulate the general conditions for character multiplicities being $\geq 2$ in terms of  $n$ and $\pi$ only.

\vspace{0.3cm}
{\large\sc Orbit Numbers:}\quad The multiplicity $c(M_{k},U)$ of the trivial $FG$-representation $U=1$ in the representation on $\CQ_{k}$ is at least as large as the number $N_{k}$ of $G$-orbits on $\CQ_{k}.$ As simple examples show, they may be different.  However, if  $p$ does not divide the order of $G$ then $M_{k}$ is semi-simple and hence $c(M_{k},1)=N_{k}.$ In this situation we may state a corollary to Theorem~\ref{Ineq2}:

\medskip
\begin{thm}\label{Ineq3} \, Let $q\geq 1$ and  suppose that $G\subseteq {\rm Aut}(\CP(n,q)).$ Suppose that the prime $p$ is co-prime to $|G|$ and $q.$  Put $\pi:=\pi(p,q).$  Then
$$
\bigl[ N_{m}\bigr]_\pi\ge
\bigl[ N_{m-1}\bigr]_\pi\ge
\bigl[ N_{m-2}\bigr]_\pi\ge
\ldots\ge
\bigl[ N_{m-s}\bigr]_\pi\ge 0\eqno (6)
$$
where $m=\lfloor\frac{n}{2}\rfloor$ and $s=\lfloor\frac{\pi}{2}\rfloor.$
\end{thm}

If one selects  $p$ large enough these inequalities include the Livingstone-Wagner theorem. However, often  improved bounds can be obtained when additional information about the prime divisors of the group is available. The following examples illustrate this. 
 
\medskip 
{\sc Example 1:\quad}
It is known from the {\sc Atlas}~\cite{atlas} that the alternating group
$A_{12}$ has a representation $G$ by $10\times 10$-matrices over ${\rm GF}(2)$.
We have $|G|=2^{9}\cdot 3^{5}\cdot 5^{2}\cdot 7\cdot 11;$
so choose $p$ to be $17, \,73$ or $127.$ It is easy to check that $\pi(17,2)=8$, 
$\pi(73,2)=9$ and $\pi(127,2)=7$, and that the corresponding system of inequalities~(6) takes the form:
$$
\left\{
\begin{array}{lcr}
N_{2}\ge N_{1}+N_{0} ,     & &(\pi=9)\medskip\\
N_3\ge N_2+N_0\ge 2N_1 ,   & &(\pi=8)\medskip\\
N_4\ge N_3+N_0\ge N_2+N_1 .& &(\pi=7)
\end{array}
\right.
$$

From these we immediately conclude that
$N_2\ge 2$, $N_3\ge 3$ and $N_4\ge 4$.
\dne

\medskip
{\sc Example 2:\quad} For $q=1$ Theorem~\ref{Ineq3} is a statement about the orbits of a permutation group $(G,\Omega)$ on the subsets of $\Omega,$ and here  $\pi(p,1)=p$. In the following example 
suppose that $n=|\Omega|=24$ and that $|G|$ is not divisible by 
$13,\, 17$ and $19.$  Then we have  the inequalities
$$
\left\{
\begin{array}{lcr}
N_{12}\ge N_{11}+N_0\ge N_{10}+N_1\ge\ldots\ge N_6+N_5\ge 2N_5 ,& \,\,(p=13)&\medskip\\
N_8\ge N_7+N_0\ge N_6+N_1\ge\ldots\ge N_4+N_3\ge 2N_3 ,         & \,\,(p=17)&\medskip\\
N_6\ge N_5+N_0\ge N_4+N_1\ge N_3+N_2\ge 2N_2 .                  & \,\,(p=19)&
\end{array}
\right.
$$
For such a group we infer that $N_7\ge N_6\ge 2$,
$N_{11}\ge\ldots\ge N_8\ge 3$ and $N_{12}\ge 4$.
For instance, the Mathieu group ${{\rm M}_{24}}$
has order $244,823,040=2^{10}\cdot 3^{3}\cdot 5\cdot 7\cdot 11\cdot23$.
Its orbit numbers are 
$N_{12}=5$, $N_{11}=\ldots=N_8=3$, $N_7=N_6=2$ and 
$N_5=\ldots=N_0=1.$ This shows that the  inequalities~(6)\, produce quite 
sharp bounds in some cases.
\dne

\end{document}